\documentclass[12pt]{article}
\usepackage{amsmath,latexsym,amssymb}
\begin{document}

\thispagestyle{empty}

\begin{center}
{\bf\large Curvature Tensor for the Spacetime of General Relativity }\\

\vspace{1cm}

{\bf Zafar Ahsan$^{*}$ and Musavvir Ali$^{**}$}

\end{center}
\vspace{1.6cm}

{\noindent\bf {Abstract :}}
  In the differential geometry of certain $F$-structures, the role of $W$-curvature tensor is very well known. A detailed study of this tensor has been made on the\linebreak spacetime of general relativity. The spacetimes satisfying Einstein field equations with vanishing $W$-tensor have been considered and the existence of Killing and\linebreak conformal Killing vector fields has been established. Perfect fluid spacetimes with \linebreak vanishing $W$-tensor have also been considered. The divergence of $W$-tensor is studied in\linebreak detail and it is seen, among other results, that a perfect fluid spacetime  with\linebreak conserved $W$-tensor represents either an Einstein space or a Friedmann-Robertson-Walker cosmological model. \\


\

\vspace{1.6cm}

{\noindent {$*$ Permanent Address: Department of Mathematics, Aligarh Muslim University, Aligarh-$~~~~$202002 (India).\newline $~~~~$e-mail: zafar.ahsan@rediffmail.com\\
\newline$**$Corresponding Author: Department of Mathematics, Aligarh Muslim University, Aligarh-202002 (India).\newline $~~~~$e-mail: musavvir.alig@gmail.com}}\\

\newpage
{\noindent\bf {1. Introduction }}\\

Pokhariyal and Mishra [8] have introduced a new curvature tensor and studied its properties. This (0,4) type tensor, denoted by $W_2$, is defined as

$$\begin{array}{lll}
W_2(X,Y,Z,T)=R(X,Y,Z,T)+\frac{1}{n-1}[g(X,Z)\mbox{Ric}(Y,T)-g(Y,Z)\mbox{Ric}(X,T)]\end{array}\leqno(1.1)$$

{\noindent {where $$R(X,Y,Z,T)=g(R(X,Y,Z),T)~~\mbox{and}~~R(X,Y,Z)=D_XD_YZ-D_YD_XZ-D_{[X,Y]}Z$$ ($D$ being the Riemannian connection) is the Riemann curvature tensor,\linebreak $\mbox{Ric}(X,Y)=g(R(X),Y)$ is the (0,2) type Ricci tensor and $R$ is the scalar\linebreak curvature.\\

For the sake of convenience, we shall denote this tensor by W. This curvature tensor, in local coordinates, can be expressed as
$$\begin{array}{lll}
W_{abcd}=R_{abcd}+\frac{1}{n-1}[g_{ac}R_{bd}-g_{bc}R_{ad}]\end{array}\leqno(1.1a)$$

\noindent and satisfies the following properties:
$$\begin{array}{lll}
W_{abcd}=-W_{bacd}, ~W_{abcd}\not=-W_{abdc},~ W_{abcd}\not=W_{cdab}\end{array}\leqno(1.2)$$
$$\begin{array}{lll}
W_{abcd}+W_{bcad}+W_{cabd}=0\end{array}\leqno(1.3)$$

\noindent which, in index-free notation, can be expressed as
$$\begin{array}{lll}
&W(X,Y,Z,T)=-W(Y,X,Z,T)\\
&W(X,Y,Z,T)\not=-W(X,Y,T,Z)\\
&W(X,Y,Z,T)\not=W(Z,T,X,Y)\end{array}\leqno(1.2a)$$

$$\begin{array}{lll}
W(X,Y,Z,T)+W(Y,Z,X,T)+W(Z,X,Y,T)=0\end{array}\leqno(1.3a)$$

In the differential geometry of certain $F$-structures, this tensor ($W_2$ or $W$) has extensively been studied by a number of workers. Thus for example, for a Sasakian manifold this tensor was studied by Pokhriyal [9]; while for a P-Sasakian manifold\linebreak Matsumoto et al [7] have studied this tensor. On the other hand, in terms of\linebreak $W_2$-tensor, Shaikh et al [10] have introduced the notion of weekly $W_2$-symmetric \linebreak manifolds and studied their properties along with several non-trivial examples. The role of $W_2$-tensor in the study of Kenmotsu manifolds has been investigated by Yildiz and De [14] while N(k)-quasi Einstein manifolds satisfying the conditions $R(\xi, X).W_2=0$ have been considered by Taleshian and Hosseinzadeh [12]. Most recently, Venkatesha et al [13] have studied Lorentzian para-Sasakian manifolds satisfying certain conditions on $W_2$-curvature tensor.\\

Motivated by the all important role of $W_2$-curvature tensor in the study of\linebreak certain differential geometric structures, in this paper we have made a detailed study of this tensor on the spacetime of general relativity. In section 2, algebraic properties of\linebreak $W$-curvature tensor are given, while the spacetimes with vanishing $W$-tensor have been considered in section 3. The existence of Killing and conformal Killing vector fields has been established for such spacetimes. Perfect fluid spacetimes satisfying Einstein field equations for $W$-flat spaces have also been studied. Section 4 deals with a detailed study of divergence of $W$-curvature tensor and perfect fluid spacetimes. A number of results concerning the vanishing of the divergence of $W$-tensor have been proved and it is seen that a perfect fluid spacetime with vanishing divergence of $W$-tensor is\linebreak either an Einstein space or a Friedmann-Robertson-Walker cosmological model. \\

{\noindent\bf {2. $W$-Curvature Tensor}}\\

It is known that the Bianchi differential identity is given by
$$\begin{array}{lll}
\nabla_eR_{abcd}+\nabla_cR_{abde}+\nabla_dR_{abec}=0\end{array}\leqno(2.1)$$

\noindent where a semi-colon denotes the covariant differentiation. This equation, in index-free notation can be expressed as
$$\begin{array}{lll}
(\nabla_UR)(X,Y,Z,T)+(\nabla_ZR)(X,Y,T,U)+(\nabla_TR)(X,Y,U,Z)=0\end{array}\leqno(2.1a)$$

Let $V_4$ be the 4-dimensional spacetime of general relativity, then equation (1.1) takes the form
$$\begin{array}{lll}
W(X,Y,Z,T)=R(X,Y,Z,T)+\frac{1}{3}[g(X,Z)\mbox{Ric}(Y,T)-g(Y,Z)\mbox{Ric}(X,T)]\end{array}\leqno(2.2)$$

From equations (2.1a) and (2.2a) it is easy to see that $W$-curvature tensor satisfies the equation

$$\begin{array}{lll}
(\nabla_XW)(Y,Z,T,U)&+~(\nabla_YW)(Z,X,T,U)~+~(\nabla_ZW)(X,Y,T,U)\\
&=\frac{1}{3}[g(Y,T)(\nabla_X\mbox{Ric})(Z,U)-g(Z,T)(\nabla_X\mbox{Ric})(Y,U)\\
&+g(Z,T)(\nabla_Y\mbox{Ric})(X,U)-g(X,T)(\nabla_Y\mbox{Ric})(Z,U)\\
&+g(X,T)(\nabla_Z\mbox{Ric})(Y,U)-g(Y,T)(\nabla_Z\mbox{Ric})(X,U)]\end{array}\leqno(2.3)$$

\noindent which, in local coordinates, can be expressed as
$$\begin{array}{lll}
&\nabla_aW_{bcde}+\nabla_bW_{cade}+\nabla_cW_{abde}\\
&=\frac{1}{3}[g_{bd}(\nabla_aR_{ce}-\nabla_cR_{ae})+g_{cd}(\nabla_bR_{ae}-\nabla_aR_{be})+g_{ad}(\nabla_cR_{be}-\nabla_bR_{ce})]\end{array}\leqno(2.4)$$

If the Ricci tensor $R_{ab}$ is of Codazzi type [4], then
$$\begin{array}{lll}
(\nabla_X\mbox{Ric})(Y,Z)=(\nabla_Y\mbox{Ric})(X,Z)=(\nabla_Z\mbox{Ric})(X,Y)\end{array}\leqno(2.5)$$

\noindent or, in local coordinates
$$\begin{array}{lll}
&\nabla_aR_{ce}=\nabla_cR_{ae}=\nabla_eR_{ac}\end{array}\leqno(2.5a)$$

\noindent(The geometrical and topological consequences of the existence of a non-trivial Codazzi tensor on a Riemannian manifold have been given in [4].The simplest Codazzi tensors are parallel one.)\\

From equation (2.5a), equation (2.4) leads to
$$\begin{array}{lll}
\nabla_aW_{bcde}+\nabla_bW_{cade}+\nabla_cW_{abde}=0\end{array}\leqno(2.6)$$

\noindent which we call as Bianchi-like identity for $W$-curvature tensor.\\

\noindent Conversely, if $W$-curvature tensor satisfies the Bianchi-like identity (2.6), then equation (2.4) reduces to
$$g_{bd}(\nabla_aR_{ce}-\nabla_cR_{ae})+g_{cd}(\nabla_bR_{ae}-\nabla_aR_{be})+g_{ad}(\nabla_cR_{be}-\nabla_bR_{ce})=0$$

\noindent Contraction of the equation with $g^{be}$ yields
$$\nabla_aR_{dc}=\nabla_cR_{da}$$
\noindent which shows that the Ricci tensor is Codazzi. Thus we have\\

{\noindent\bf {Theorem 2.1 :}} For a $V_4$, the Ricci tensor is of Codazzi type if and only if $W$-curvature tensor satisfies the Bianchi-like identity (2.5).\\

{\noindent\bf{3. Spacetimes with vanishing $W$-curvature Tensor}}\\

Let $V_4$ be the spacetime of general relativity, then from equation (1.1a) we have
$$W_{abcd}=R_{abcd}+\frac{1}{3}[g_{ac}R_{bd}-g_{bc}R_{ad}]\leqno(3.1)$$
$$W^h_{bcd}=R^h_{bcd}+\frac{1}{3}[{\delta}^h_cR_{bd}-g_{bc}R^h_d]\leqno(3.2)$$

\noindent Contraction of equation (3.2) over $h$ and $d$ leads to
$$W_{bc}=\frac{4}{3}(R_{bc}-\frac{1}{4}Rg_{bc})\leqno(3.3)$$

\noindent so that
$$g^{bc}W_{bc}=0\leqno(3.4)$$
{\noindent\bf {Definition 3.1 :}} A spacetime is said to be $W$-flat if its $W$-curvature tensor vanishes.\\

For a $W$-flat spacetime, equation (3.2) leads to
$$R^h_{bcd}=-\frac{1}{3}(\delta^h_cR_{bd}-g_{bc}R^h_d)\leqno(3.5)$$
\noindent which on contraction over $h$ and $d$ yields
$$R_{bc}=\frac{R}{4}g_{bc}\leqno(3.6)$$
\noindent Equation (3.6) shows that a $W$-flat spacetime is an Einstein space. Thus, we have\\

{\noindent\bf{Theorem 3.1 :}} A $W$-flat spacetime is an Einstein space and consequently the scalar curvature $R$ is a covariantly constant, i.e. $\nabla_lR=0$.\\

{\noindent\bf{Remark 3.1 :}} The spaces of constant curvature play a significant role in cosmology. The simplest cosmological model is obtained by assuming that the universe is isotropic and homogenous. This is known as cosmological principle. This principle, when\linebreak translated into the language of Riemannian geometry, asserts that the three\linebreak dimensional position space is a space of maximal symmetry [2, 11], that is space of constant curvature whose curvature depends upon time. The cosmological solutions of Einstein equations which contain a three dimensional space-like surface of a constant curvature are the Robertson-Walker metrics, while four dimensional space of constant curvature is the de Sitter model of the universe (c.f., [11]).\\

In general theory of relativity, the curvature tensor describing the gravitational field consists of two parts viz., the matter part and the free gravitational part. The interaction between these parts is described through Bianchi identities. For a given distribution of matter, the construction of gravitational potential satisfying Einstein field equations is the principal aim of all investigations in gravitational physics; and this has often been achieved by imposing symmetries on the geometry compatible with the dynamics of the chosen distribution of matter. The geometrical symmetries of spacetime are expressed  through the equation
$$\pounds_{\xi}A-2\Omega A=0\leqno(3.7)$$
\noindent where $A$ represents a geometrical/physical quantity, $\pounds_{\xi}$ denotes the Lie derivative with respect to a vector field $\xi$ and $\Omega$ is a scalar.

Let the Einstein field equations with a cosmological term be
$$R_{bc}-\frac{1}{2}Rg_{bc}+\Lambda g_{bc}=\mbox{k} T_{bc}\leqno(3.8)$$

\noindent where $R_{bc}$ is the Ricci tensor, $R$, the scalar curvature, $\Lambda$, the cosmological constant, k, the non-zero gravitational constant and $T_{bc}$, the energy-momentum tensor. From equation (3.6), equation (3.8) leads to
$$(\Lambda-\frac{R}{4})g_{bc}=\mbox{k} T_{bc}\leqno(3.9)$$

\noindent Since for $W$-flat spacetime $R$ is constant, taking the Lie derivative of both sides of equation (3.9) we get

$$(\Lambda-\frac{R}{4})\pounds_{\xi}g_{bc}=\mbox{k}\pounds_{\xi}T_{bc}\leqno(3.10)$$
We thus have the following\\

{\noindent\bf{Theorem 3.2 :}} For a $W$-flat spacetime satisfying Einstein field equations with a \linebreak cosmological term, there exists a Killing vector field $\xi$ if and only if the Lie derivative of the energy momentum tensor vanishes with respect to $\xi$.\\

{\noindent\bf{Definition 3.2 :}}  A vector field $\xi$ satisfying the equation
$$\pounds_\xi g_{bc}=2\Omega g_{bc}\leqno(3.11)$$
\noindent is called a conformal Killing vector field, where $\Omega$ is a scalar. A spacetime satisfying equation (3.11) is said to admit a conformal motion.

From equation (3.10) and (3.11) we have
$$2\Omega (\Lambda-\frac{R}{4})g_{bc}=\mbox{k}\pounds_\xi T_{bc}\leqno(3.12)$$

\noindent which on using (3.9) leads to
$$\pounds_\xi T_{bc}=2\Omega T_{bc}\leqno(3.13)$$

The energy-momentum tensor $T_{bc}$ satisfying equation (3.13) is said to posses the symmetry inheritance property (cf; [1]). Thus, we can state the following\\

{\noindent\bf{Theorem 3.3 :}}  A $W$-flat spacetime satisfying the Einstein field equations with a\linebreak cosmological term admits a conformal Killing vector field if and only if the\linebreak energy-momentum tensor has the symmetry inheritance property.\\

Consider now a perfect fluid spacetime with vanishing $W$-curvature tensor. The energy-momentum tensor $T_{bc}$ for a perfect fluid is given by
$$T_{bc}=(\mu+p)u_bu_c+pg_{bc}\leqno(3.14)$$
\noindent where $\mu$ is the energy density, $p$ the isotropic pressure, $u_a$ the velocity of fluid such that $u_au^a=-1$ and $g_{bc}u^b=u_c$.

From equations (3.9) and (3.14) we have
$$(\Lambda-\frac{R}{4}-\mbox{k} p)g_{bc}=\mbox{k}(\mu+p)u_bu_c\leqno(3.15)$$
\noindent which on multiplication  with $g^{bc}$ yields
$$R=\mbox{k}(\mu-3p)+4\Lambda\leqno(3.16)$$
\noindent Also, the contraction of equation (3.15) with $u^bu^c$ leads to
$$R=4(\mbox{k}\mu+\Lambda)\leqno(3.17)$$
\noindent A comparison of equations (3.16) and (3.17) now yields
$$\mu+p=0\leqno(3.18)$$
\noindent which means that either $\mu=0$, $p=0$ (empty spacetime) or the perfect fluid spacetime satisfies the vacuum-like equation of state [6]. We thus have the following\\

{\noindent\bf{Theorem 3.4 :}} For a $W$-flat perfect fluid spacetime satisfying Einstein field equations with a cosmological term, the matter contents of the   spacetime obey the vacuum-like equation of state.\\

The Einstein field equations in the presence of matter are given by
$$R_{ab}-\frac{1}{2}Rg_{ab}=\mbox{k} T_{ab}\leqno(3.19)$$
\noindent which on multiplication with $g^{ab}$ leads to
$$R=-\mbox{k} T\leqno(3.20)$$
\noindent Equation (3.19) may be expressed as
$$R_{ab}=\mbox{k}(T_{ab}-\frac{1}{2}Tg_{ab})\leqno(3.21)$$
\noindent so that in the absence of matter, the field equations are
$$R_{ab}=0\leqno(3.22)$$
\noindent Equations (3.22) are the field equations for empty spacetime.\\

It is known that [2, 11] the energy-momentum tensor for the electromagnetic field is given by
$$T_{ab}=-F_{ac}F^c_b+\frac{1}{4}g_{ab}F_{pq}F^{pq}\leqno(3.23)$$
\noindent Where $F_{ac}$ represents the skew-symmetric electromagnetic field tensor satisfying Maxwell's equations. From equation (3.23) it is evident that $T^a_a=T=0$ and the Einstein\linebreak equations for a purely electromagnetic distribution take the form
$$R_{ab}=\mbox{k} T_{ab}\leqno(3.24)$$
\noindent Moreover, from equation (3.19) we have
$$\nabla_cR_{ab}-\nabla_bR_{ac}=\mbox{k} (\nabla_cT_{ab}-\nabla_bT_{ac})+\frac{1}{2}(g_{ab}\nabla_cR-g_{ac}\nabla_bR)\leqno(3.25)$$
\noindent If $T_{ab}$ is of Codazzi type, then equation (3.25) becomes
$$\nabla_cR_{ab}-\nabla_bR_{ac}=\frac{1}{2}(g_{ab}\nabla_cR-g_{ac}\nabla_bR)\leqno(3.26)$$

\noindent which on multiplication with $g^{ab}$, after simplification, leads to
$$\nabla_bR^{ab}=0\leqno(3.27)$$
\noindent Thus we have\\

{\noindent\bf{Theorem 3.5 :}}  For a $V_4$ satisfying Einstein field equations, the Ricci tensor is\linebreak conserved if the energy-momentum tensor is Codazzi.\\

{\noindent\bf{4. Divergence of $W$-curvature tensor and perfect fluid spacetimes}}\\

The Bianchi identities are given by
$$\nabla_eR^h_{bcd}+\nabla_cR^h_{bde}+\nabla_dR^h_{bec}=0\leqno(4.1)$$
\noindent Contracting equation (4.1) over h and e, using the symmetry properties of Riemann curvature tensor, we get
$$\nabla_hR^h_{bcd}=\nabla_dR_{bc}-\nabla_cR_{bd}\leqno(4.2)$$

It is known that Riemannian manifolds for which the divergence of the\linebreak curvature tensor vanish identically are known as manifolds with harmonic curvature. The curvature of such manifolds occur as a special case of Yang-Mills fields. These manifolds also form a natural generalization of Einstein spaces and of conformally flat manifolds with constant scalar curvature [3]. From equation (4.2) we thus have\\

{\noindent\bf{Theorem 4.1 :}}  If the Ricci tensor is of Codazzi type then the manifold $V_4$ \linebreak is of harmonic curvature and conversely.\\

{\noindent\bf{Remark 4.1 :}}  The Ricci tensor $\mbox{Ric}(X,Y)$ is said to be parallel if
$$\nabla_Z\mbox{Ric}(X,Y)-\nabla_Y\mbox{Ric}(X,Z)=0$$
\noindent which means that the simplest Codazzi tensors are parallel ones.\\

Now from equation (3.2) we have
$$\nabla_eW^h_{bcd}=\nabla_eR^h_{bcd}+\frac{1}{3}(\delta^h_c\nabla_eR_{bd}-g_{bc}\nabla_eR^h_{d})\leqno(4.3)$$
\noindent so that the divergence of $W$-curvature tensor is given by

$$\nabla_h W^h_{bcd}=\nabla_h R^h_{bcd}+\frac{1}{3}(\nabla_c R_{bd}-g_{bc}\nabla_h R^h_d)\leqno(4.4)$$
\noindent which leads to the following\\

{\noindent\bf{Theorem 4.2 :}} For a spacetime possessing harmonic curvature with divergence-free $W$-tensor, the Ricci tensor is covariantly constant.\\

While from equations (4.2) and (4.4) we have
$$\nabla_hW^h_{bcd}=\nabla_dR_{bc}-\nabla_cR_{bd}+\frac{1}{3}(\nabla_cR_{bd}-g_{bc}\nabla_hR^h_d)\leqno(4.5)$$
\noindent Thus we can state the following\\

{\noindent\bf{Theorem 4.3 :}}  If for a $V_4$ the divergence of $W$-tensor vanishes and the Ricci tensor is covariantly constant then the spacetime is of constant curvature.\\

Now using equation (3.21) in equation (4.5) we get
$$\nabla_hW^h_{bcd}=\mbox{k}(\nabla_dT_{bc}-\frac{2}{3}\nabla_cT_{bd})+\frac{\mbox{k}}{3}(g_{bd}\nabla_cT-\frac{5}{2}g_{bc}\nabla_dT)\leqno(4.6)$$
\noindent so that for a purely electromagnetic distribution, we have
$$\nabla_hW^h_{bcd}=\mbox{k}(\nabla_dT_{bc}-\frac{2}{3}\nabla_cT_{bd})\leqno(4.7)$$
\noindent which leads to\\

{\noindent\bf{Theorem 4.4 :}} For a spacetime satisfying the Einstein equations for a purely\linebreak electromagnetic distribution, the $W$-curvature tensor is conserved if the energy-momentum tensor is covariantly constant and conversely.\\

Consider now the energy-momentum tensor for a perfect fluid [cf, equation (3.14)]
$$T_{ab}=(\mu+p)u_au_b+pg_{ab}\leqno(4.8)$$
\noindent which leads to
$$T=-\mu+3p\leqno(4.9)$$
 If $T_{ab}$ is Codazzi and $\nabla_hW^h_{bcd}=0$ then from equations (4.8) and (4.9), equation (4.6) leads to $(\mbox{k}=1)$
$$\begin{array}{lll}
&\frac{1}{3}[\nabla_c(\mu+p)u_bu_d+(\mu+p)\nabla_cu_{b}u_d+(\mu+p)u_b\nabla_cu_{d}+\nabla_cpg_{bd}]\\
&+\frac{1}{3}g_{bd}\nabla_c(-\mu+3p)-\frac{5}{6}g_{bc}\nabla_d(-\mu+3p)=0\end{array}\leqno(4.10)$$
\noindent Since $\nabla_bu_{a}u^a=0$, contracting equation (4.10) with $g^{bd}$ we get
$$\nabla_c(\mu-3p)=0\leqno(4.11)$$
\noindent that is
$$(\mu-3p)=\mbox{constant}\leqno(4.12)$$
\noindent Thus, we have\\

{\noindent\bf{Theorem 4.5 :}} If for a perfect fluid spacetime, the divergence of $W$-curvature tensor vanishes and the energy-momentum tensor is of Codazzi type $(\mu-3p)$ is constant.\\

It is known that [5] for a radiative perfect fluid spacetime $(\mu=3p)$ the resulting universe is isotropic and homogenous. Thus by choosing the constant in equation (4.12) as zero, we have the following\\

{\noindent\bf{Corollary 4.1 :}}  If the energy-momentum tensor for a divergence-free $W$-fluid\linebreak spacetime is of Codazzi type then the resulting spacetime is radiative and consequently isotropic and homogenous.\\

Now consider the spacetime for which the divergence of $W$-curvature tensor\linebreak vanishes, then from equation (4.6), we have $(\mbox{k}=1)$
$$\nabla_dT_{bc}-\frac{5}{6}g_{bc}\nabla_dT=\frac{2}{3}\nabla_cT_{bd}-\frac{1}{3}g_{bd}\nabla_cT\leqno(4.13)$$
\noindent which on using equation (4.8) and (4.9) leads to

$$\begin{array}{lll}
&\nabla_d(\mu+p)u_bu_c+(\mu+p)\nabla_du_{b}u_c+(\mu+p)u_b\nabla_du_{c}+\nabla_dpg_{bc}-\frac{5}{6}g_{bc}\nabla_d(-\mu+3p)\\
&=\frac{2}{3}[\nabla_c(\mu+p)u_bu_d+(\mu+p)\nabla_cu_{b}u_d+(\mu+p)u_b\nabla_cu_{d}+\nabla_cpg_{bd}]-\frac{1}{3}g_{bd}\nabla_c(-\mu+3p)\end{array}$$
Contracting this equation with $u^d$, we get
$$\begin{array}{lll}
&(\mu+p\dot{)}u_bu_c+(\mu+p)\dot{u}_bu_c+(\mu+p)u_b\dot{u}_c\\
&+\dot{p}g_{bc}-\frac{5}{6}(-\mu+3p\dot{)}g_{bc}+\frac{2}{3}\nabla_c(\mu+p)u_b+\frac{2}{3}(\mu+p)\nabla_cu_{b}\\
&-\frac{2}{3}\nabla_cpu_b+\frac{1}{3}\nabla_c(-\mu+3p)u_b=0\end{array}\leqno(4.14)$$
\noindent where an over head dot denotes the covariant derivative along the fluid flow vector $u_a$ (that is, $(\mu+p\dot{)}=\nabla_c(\mu+p)u^c$, $\dot{u}_b=\nabla_cu_{b}u^c$, $\dot{p}=\nabla_dpu^d$, $\nabla_bu_{a}u^a=0$, etc.).\\

Also, the conservation of energy-momentum tensor $(\nabla_bT^{ab}=0)$ leads to
$$(\mu+p)\dot{u}_a=-\nabla_ap+\dot{p}u_a~~~~~~~~~~~~ \mbox{(force~equation)}\leqno(4.15)$$
$$\dot{\mu}=-(\mu+p)\nabla_au^a=-(\mu+p)\theta~~~~~~\mbox{(energy~equation)}\leqno(4.16)$$
\noindent Moreover, the covariant derivative of the velocity vector can be splitted into kinematical quantities as
$$\nabla_bu_{a}=\frac{1}{3}\theta(g_{ab}+u_au_b)-\dot{u}_au_b+\sigma_{ab}+\omega_{ab}\leqno(4.17)$$
\noindent where $\theta=\nabla_au^a$, is the expansion scalar, $\dot{u}_a=\nabla_bu_{a}u^b$, the acceleration vector $\sigma_{ab}=h^c_ah^d_bu_{(c;d)}-\frac{1}{3}\theta h_{ab}$, the symmetric shear tensor $(h_{ab}=g_{ab}-u_au_b)$ and $\omega_{ab}=h^c_ah^d_bu_{[c;d]}$ is the skew symmetric vorticity or rotation tensor.\\

Using force equation (4.15) in equation (4.14), we get
$$\begin{array}{lll}
(\mu-p\dot{)}u_bu_c-\nabla_bpu_c+\dot{p}g_{bc}+\frac{2}{3}(\mu+p)\nabla_cu_b-\frac{5}{6}(-\mu+3p\dot{)}g_{bc}+\frac{1}{3}\nabla_c\mu u_b=0\end{array}\leqno(4.18)$$

\noindent Contracting this equation with $u^b$, we get
$$\frac{1}{6}(\mu-3p\dot{)}u_c-\frac{1}{3}\nabla_c\mu=0\leqno(4.19)$$
Thus we have the following\\

{\noindent\bf{Theorem 4.6 :}}  For a perfect fluid spacetime with divergence-free $W$-curvature tensor, the pressure and density of the fluid are constant.\\

Now contracting the equation (4.18) with $u^c$, we get
$$\frac{3}{2}(\mu-3p\dot{)}u_b+\frac{2}{3}(\mu+p)\dot{u}_b+\nabla_bp=0\leqno(4.20)$$
\noindent which on using force equation (4.15) leads to
$$\frac {3}{2}\dot{\mu}u_b-\frac{31}{6}\nabla_dpu^du_b+\frac{1}{3}\nabla_bp=0$$
From energy equation (4.16), this equation yields
$$\frac{3}{2}[(\mu+p)\theta]u_b+\frac{31}{6}\nabla_dpu^du_b-\frac{2}{3}\nabla_bp=0\leqno(4.21)$$

We thus have the following\\

{\noindent\bf{Theorem 4.7 :}}  For a perfect fluid spacetime having conserved $W$-curvature tensor, either pressure and density of the fluid are constant over the space-like hypersurface orthogonal to the fluid four velocity or the fluid is expansion-free.\\

Consider now a perfect fluid spacetime for which $\nabla_hW^h_{bcd}=0$, then from equations (4.14) and (4.15), we have
$$\begin{array}{lll}
&(\mu+3p\dot{)}u_bu_c+\dot{p}g_{bc}-\nabla_bpu_c-\frac{5}{6}(-\mu+3p\dot{)}g_{bc}\\
&+\frac{1}{3}\nabla_c(2\mu-3p\dot{)}u_b+\frac{2}{3}(\mu+p)\nabla_cu_{b}+\frac{1}{3}\nabla_c(-\mu+3p)u_b=0\end{array}\leqno(4.22)$$
\noindent which on contraction with $u^b$ leads to
$$\begin{array}{lll}
(\mu+3p\dot{)}u_c+\frac{5}{6}(-\mu+3p\dot{)}u_c+\frac{1}{3}\nabla_c(2\mu-3p\dot{)}+\frac{1}{3}\nabla_c(-\mu+3p)=0\end{array}\leqno(4.23)$$

\noindent This equation is satisfied only when
$$\nabla_c(2\mu-3p)=-3(\mu+3p\dot{)}u_c\leqno(4.24)$$
\noindent and
$$\nabla_c(-\mu+3p)=-\frac{5}{2}(-\mu+3p\dot{)}u_c\leqno(4.25)$$

Using force equation (4.15) in equation (4.24), we get
$$(\mu+p)\dot{u}_c=-\frac{1}{3}\dot{\mu}u_c-\frac{2}{9}\nabla_c(\mu+3p)\leqno(4.26)$$
\noindent Also, from force equation (4.15) and equation (4.25), we have
$$\nabla_c\mu=-3(\mu+p)\dot{u}_c-\frac{1}{2}(5\mu-21p\dot{)}u_c\leqno(4.27)$$

From equations (4.22) and (4.24) we have
$$\begin{array}{lll}
&\dot{p}g_{bc}-\nabla_bpu_c-\frac{5}{6}(-\mu+3p\dot{)}g_{bc}\\
&+\frac{1}{3}(\mu+p)\nabla_cu_{b}+\frac{1}{3}\nabla_c(-\mu+3p)u_b=0\end{array}\leqno(4.28)$$
\noindent which on multiplication with $g^{bc}$ leads to
$$\dot{p}=\frac{10}{9}(-\mu+3p\dot{)}-\frac{1}{9}(\mu+p)\theta-\frac{1}{9}\nabla_c(-\mu+3p)u^c\leqno(4.29)$$

Using energy equation (4.16) this equation yields
$$\dot{p}=\frac{3}{7}\dot{\mu}+\frac{1}{21}\nabla_c(-\mu+3p)u^c\leqno(4.30)$$

Also, from the equations (4.15) and (4.28), we have
$$\begin{array}{lll}
&(\mu+p)\dot{u}_bu_c+\dot{p}(g_{bc}+u_bu_c)+\frac{1}{3}(\mu+p)\nabla_cu_{b}\\
&-\frac{5}{6}(-\mu+3p\dot{)}g_{bc}+\frac{1}{3}\nabla_c(-\mu+3p)u_b=0\end{array}\leqno(4.31)$$
\noindent which on using equations (4.25) and (4.29) yields
$$(\mu+p)[3\dot{u}_bu_c-\frac{1}{3}\theta(g_{bc}+u_bu_c)+\nabla_cu_{b}]=0\leqno(4.32)$$

This equation suggests that either
$$\mu+p=0\leqno(4.33)$$
\noindent or
$$\nabla_cu_{b}+3\dot{u}_bu_c-\frac{1}{3}\theta(g_{bc}+u_bu_c)=0\leqno(4.34)$$

From equation (4.33) either $\mu=0$, $p=0$ (neither matter nor radiation) or the perfect fluid spacetime with $\nabla_hW^h_{bcd}=0$ satisfies the vacuum-like equation of state [6].\\

Moreover, from equations (4.16), (4.25) and (4.29) we have
$$-\mu+3p=\mbox{constant}\leqno(4.35)$$
\noindent while from equations (4.17) and (4.34) we have
$$2\dot{u}_bu_c+\sigma_{bc}+\omega_{bc}=0\leqno(4.36)$$
\noindent which is satisfied only when $\dot{u}_b=0$, $\sigma_{bc}=0$, $\omega_{bc}=0$\\

Therefore, if $\mu+p\not=0$ then the above discussions show that the fluid is shear-free, rotation-free, acceleration-free and the energy density and pressure are constants over the space-like hypersurface orthogonal to the fluid four velocity.

It may be noted that vanishing of acceleration, shear and rotation and the\linebreak constantness of energy density and pressure over the space-like hypersurface\linebreak orthogonal to the fluid flow vector are the conditions for a spacetime to represent Friedmann-Robertson-Walker cosmological model provided that $\mu+p\not=0$\\

Hence, summing up these discussions, we can state the following\\

{\noindent\bf{Theorem 4.8 :}} A perfect fluid spacetime with conserved $W$-curvature tensor is\linebreak either an Einstein spacetime $(\mu+p)=0$ or a Friedmann-Robertson-Walker\linebreak cosmological model satisfying $\mu-3p=\mbox{constant}$.\\

{\noindent\bf\large\underline{Acknowledgment}}

The research of M. Ali was supported by UGC Startup grant No.F.30-90/2015(BSR).

\vspace{.5cm}

{\noindent\bf\large\underline{References}}

\begin{enumerate}
\item[${[1]}$]   Ahsan, Z. : On a geometrical symmetry of the spacetime of general relativity. {\it Bull. Cal. Math. Soc.}, 97 (3) 191-200, 2005.
\item[${[2]}$]   Ahsan, Z. : Tensor Analysis with Application. {\it Anshan Ltd, Tunbridge Wells, United Kingdom,} 2008.
\item[${[3]}$]   Besse, A.L. : Einstein Manifolds. {\it Springer} (1987) 235.
\item[${[4]}$]   Derdzinski, A. and Shen, C.L. : Codazzi tensor fields, Curvature and Pontryagin form. {\it Proc. London Math. Society}, 47 (3) 15-26,  1983.
\item[${[5]}$]   Ellis, G.F.R. : Relativistic Cosmology in General Relativity and Cosmology. Edited by R.K. Sachs. {\it Academic Press (London)}104-179, 1971.
\item[${[6]}$]   Kallingas, D., Wesson, P. and Everitt, C.W.F. : Flat FRW models with variable $G$ and $\Lambda$. {\it Gen. Relativity Gravit.} 24(4), 351-357, 1992.
\item[${[7]}$]   Matsumoto, K., Ianus, S. and Milrai, I. : On P-Sasakian manifolds which admit certain tensor fields. {\it Publicationes Mathematicae Debrecen,} 33(3-4)199-204, 1986.
\item[${[8]}$]   Pokhriyal, G.P. and Mishra, R.S. : Curvature Tensors and their relativistic significance. {\it Yokohama Math. Journal} 18,105-108, 1970.
\item[${[9]}$]   Pokhriyal, G.P. : Study of a new curvature tensor in a Sasakian manifold.{\it Tensor,} 36(2) 222-226, 1982.
\item[${[10]}$]   Shaikh, A.A., Jana, S.K. and Eyasmin, S. : On weakly $W_2$-symmetric manifolds. {\it Sarajevo Journal of Mathematics,} 3(15)73-91, 2007.
\item[${[11]}$]  Stephani, H. : General Relativity- An Introduction to the Theory of Gravitational Field. {\it Cambridge University Press, Cambridge,}  1982.
\item[${[12]}$]  Teleshian, A. and Hosseinzadeh, A.A. : On $W_2$-curvature tensor N(k)-quasi\linebreak Einstein manifolds. {\it Journal of Mathematics and Computer Science,} 1(1), 28-32, 2010.
\item[${[13]}$]  Venkatesha, Bagewadi, C.S., and Pradeep Kumar, K.T. : Some Results on Lorentzian para-Sasakian manifolds. {\it ISRN Geometry 2011, Article ID 161523, 9 pages, doi: 10.5402/2011/161523.}
\item[${[14]}$]  Yildiz, A. and De, U.C. : On a type of Kenmotsu manifolds. {\it Diff. Geom-Dynamical Systems,} 12, 289-298, 2010.

\end{enumerate}
\end{document}